\documentclass[a4paper]{amsart}%
\usepackage{amsgen}
\usepackage{amsmath}
\usepackage{amstext}
\usepackage{amsbsy}
\usepackage{amsopn}
\usepackage{amsfonts}
\usepackage{amsthm}
\usepackage{amscd}
\usepackage{amsxtra}
\usepackage{eucal}
\usepackage{amscd}
\usepackage{amsfonts}
\usepackage{amssymb}
\usepackage{graphicx}%
\newtheorem{theorem}{Theorem}[section]
\newtheorem*{theorem*}{Theorem}
\newtheorem{lemma}[theorem]{Lemma}
\newtheorem{proposition}[theorem]{Proposition}
\newtheorem{corollary}[theorem]{Corollary}
\theoremstyle{definition}
\newtheorem{definition}[theorem]{Definition}

\theoremstyle{remark}
\newtheorem{remark}[theorem]{Remark}

\numberwithin{equation}{section}
\begin{document}
\title[On the cotangent cohomology of rational surface singularities]{On the cotangent cohomology of rational surface singularities with almost
reduced fundamental cycle}
\author{Trond St\o len Gustavsen}
\email{stolen@math.uio.no}

\begin{abstract}
We prove dimension formulas for the cotangent spaces $T^{1}$ and $T^{2}$ for a
class of rational surface singularities by calculating a correction term in
the general dimension formulas. We get that it is zero if the dual graph of
the rational surface singularity $X$ does not contain a particular type of
configurations, and this generalizes a result of Theo de Jong stating that the
correction term $c(X)$ is zero for rational determinantal surface
singularities. In particular our result implies that $c(X)$ is zero for
Riemenschneiders quasi-determinantal rational surface singularities, and this
also generalise results for qoutient singularities.

\end{abstract}
\maketitle

\section{Introduction}

The cotangent cohomology is important in the deformation theory of isolated
singularities. In several papers the dimensions of these modules are
calculated for classes of rational surface singularities, see for instance
\cite{bkr:inf}, \cite{bc:hypr}, \cite{js:rss}, \cite{DJ}. In \cite{cg:01}
dimension formulas for the cotangent modules $T^{1}$ and $T^{2}$ for general
rational surface singularities are given, and \cite{ste:02} gives formulas for
the higher cotangent modules. The formulas for $T^{1}$ and $T^{2}$ contains an
unavoidable correction term which vanish for large classes of rational surface
singularities (see \cite{bkr:inf}, \cite{bc:hypr}, \cite{js:rss}, \cite{DJ}),
but which seems to be difficult to compute in general. In the present paper we
investigate the correction term for rational surface singularities where the
\textit{fundamental cycle}, \cite{ar:rat}, is reduced on all non $-2$-curves.
We refer to this by saying that the fundamental cycle is \textit{almost
reduced}.

The formulas proved in \cite[Theorems 3.11 and 3.8]{cg:01} may be stated as
follows:
\begin{align*}
\dim T_{X}^{1}  &  =(e-4)+\dim T_{\widehat{X}}^{1}+c(X)\\
\dim T_{X}^{2}  &  =(e-2)(e-4)+\dim T_{\widehat{X}}^{2}+c(X)
\end{align*}
In these formulas $e$ denotes the embedding dimension of the rational surface
singularity $X$ , $\widehat{X}$ is the blowup of $X$ and $c(X)$ is the
correction term.

In section \ref{classification} we give a classification of rational surface
singularities with almost reduced fundamental cycle, in particular we define
the notion of $n$-configurations in the dual graph. The main theorem of this
paper is the following:

\begin{theorem*}
Let $X$ be a rational surface singularity with embedding dimension $e\geq4$
and with almost reduced fundamental cycle. Then the correction term $c(X)$ is
less or equal to the number of $3$-configurations in the dual graph. In
particular; $c(X)=0$ for all quasi-determinantal rational surface singularities.

Assume furthermore that the fundamental cycle intersects all non $-2$-curves
negatively. Then $c(X)$ equals the number of $3$-configurations in the dual
graph of $X.$
\end{theorem*}

Determinantal rational surface singularities have almost reduced fundamental
cycle and have only $1$-configurations in the dual graph, \cite[3.4]{wa:equ},
\cite[4.2.1]{ro:92}, \cite{DJ}, \cite{dJ:99}. In this way our result
generalize the formulas given in \cite{DJ}. Notice in particular that rational
quasi-determinantal surface singularities, as defined by Riemenschneider,
\cite{ri:81}, have almost reduced fundamental cycle, see \cite{ro:92} and
\cite{dJ:99}, and that quasi-determinantal rational surface singularities do
not have $3$-configurations in the dual graph, and thus we get that $c(X)=0$
for quasi-determinantal singularities. Remark also that the results in the
present paper, generalizes and \textit{corrects} Theorem 2.6.3 of
\cite{mythesis} where some $2$-configurations mistakenly are computed to
contribute positively to $c(X).$

\subsection*{Acknowledgment}

I thank Jan Arthur Christophersen for helpful comments. I thank RCN's
Strategic University Program in Pure Mathematics at the Dept. of Mathematics,
University of Oslo (No 154077/420) for partial financial support.

\vspace{1cm}

\section{Notation and preliminaries}

\subsection{Results and notation on rational singularities}

\label{rational}We will work over the field of complex numbers. The
singularities we study are of the form $X=\operatorname{Spec}A$ where $A=P/I$
and $P$ is a regular local $\mathbb{C}$ algebra essentially of finite type. We
will denote by ${\mathfrak{m}}$ the maximal ideal in $A={\mathcal{O}}_{X}$. A
normal surface singularity $X$ with minimal resolution $f:\widetilde
{X}\rightarrow X$ is rational if $H^{1}(\widetilde{X},{\mathcal{O}%
}_{\widetilde{X}})=0$, see \cite{ar:rat}. The exceptional divisor
$E\subset\widetilde{X}$ is a union of irreducible components $E_{i}%
\simeq{\mathbb{P}}^{1}$. There is a \textit{fundamental cycle} $Z$, supported
on $E$, defined by ${\mathfrak{m}}{\mathcal{O}}_{\widetilde{X}}$. This divisor
is the unique smallest positive divisor $Z=\sum r_{i}E_{i}$ satisfying $Z\cdot
E_{i}\leq0$ for all irreducible components $E_{i}$. The fundamental cycle may
be computed inductively as follows: Put $Z_{0}=E.$ Given $Z_{k}$ there are two possibilities:

\begin{enumerate}
\item If there is an $E_{i}$ such that $Z\cdot E_{i}>0$, then put
$Z_{k+1}=Z_{k}+E,$

\item otherwise we are finished and the fundamental cycle $Z=Z_{k}.$
\end{enumerate}

The embedding dimension of $X$, $e=\dim_{\mathbb{C}}{\mathfrak{m}%
}/{\mathfrak{m}}^{2}$, equals $-Z^{2}+1$ and the multiplicity $m(X)=e-1=-Z^{2}%
$. Note also that since $X$ is affine, then for any modification $Y\rightarrow
X$ and for any quasi-coherent sheaf $\mathcal{F}$ on $Y$ it follows from the
Theorem on Formal Functions that $H^{i}(Y,\mathcal{F})=0$ for $i\geq2,$ see
\cite[III.11.1]{AG}. We will use this fact freely throughout the paper.

The following result from \cite{tj:abs}, which shows how the blow up
$\widehat{X}$ may be obtained from $\widetilde{X}$, is important for the
considerations in this paper.

\begin{theorem}
[Tjurina]\label{Tyu} If $X$ is a rational surface singularity, then the blow
up of $X$ is isomorphic to the surface obtained from $\widetilde{X}$ by
contracting all components $E_{i}$ with $Z\cdot E_{i}=0$.
\end{theorem}

We fix the following notation. Let $\widehat{X}\overset{\pi}{\rightarrow}X$
denote the blow up of a rational surface singularity $X=\operatorname{Spec}A$
in the maximal ideal, and let $\widetilde{X}\overset{\widetilde{\pi}%
}{\rightarrow}\widehat{X}$ be the minimal resolution. Denote by $C\subset
\widehat{X}$ the exceptional divisor defined by ${\mathfrak{m}}{\mathcal{O}%
}_{\widehat{X}},$ and by $Z\subset\widetilde{X}$ the fundamental cycle. Let
$\Theta_{\widehat{X}}$ be the tangent sheaf on $\widehat{X}$, and let
$\Theta_{\widetilde{X}}$ be the tangent sheaf on $\widetilde{X}$.

\subsection{Cotangent cohomology}

\label{cotangent} We refer to \cite{an:hom} for precise definitions, but see
also \cite{cg:01} for a review of the basic facts and notation.

Assume that $S$ is a noetherian ring and that $A$ is an $S$ algebra
essentially of finite type. For an $A$ module $M$ we will consider the
\textit{cotangent cohomology} modules $T^{i}(A/S;M)=:T_{A}^{i}(M)$ and
$T_{A}^{i}(A)=:T_{A}^{i}=:T_{X}^{i}$ if $X=\operatorname{Spec}A$. If
$\mathcal{S}$ is a sheaf of rings on a scheme $Y$, $\mathcal{A}$ is an
$\mathcal{S}$ algebra and $\mathcal{M}$ is an $\mathcal{A}$-module, there are
cotangent cohomology sheaves ${\mathcal{T}}_{{\mathcal{A}}/{\mathcal{S}}}%
^{i}(\mathcal{M})$ and cotangent cohomology groups $T_{{\mathcal{A}%
}/{\mathcal{S}}}^{i}(\mathcal{M}).$ The sheaf ${\mathcal{T}}_{{\mathcal{A}%
}/{\mathcal{S}}}^{i}(\mathcal{M})$ is locally for an open affine $U\subset Y$
given as $T^{i}({\mathcal{A}}(U)/{\mathcal{S}}(U);\mathcal{M}(U)).$

For the purpose of calculations we note that if $P$ is a polynomial
$S$-algebra (or a the localization of such) mapping onto $A$ so that $A\cong
P/I$ for an ideal $I$, then $T^{1}(A/S;M)$ is the cokernel of the natural map
$\operatorname*{Der}_{S}(P,M)\rightarrow\operatorname*{Hom}_{P}(I,M).$

\subsection{The sheaves $\mathcal{F}^{i}$ and the correction term $c(X)$.}

\label{Computation}To shorten notation we define the sheaves ${\mathcal{F}%
}^{i}:={\mathcal{T}}_{\pi^{-1}A}^{i}({\mathcal{O}}_{\widehat{X}})\,$on
$\widehat{X}$.

\begin{definition}
The correction term $c(X)$ is defined as $c(X)=h^{1}(\widehat{X}%
,\mathfrak{m}\mathcal{F}^{1}).$
\end{definition}

The importance of this invariant is given by the following formulas which are
valid for all rational surface singularities, see \cite[Theorems 3.11 and
3.8]{cg:01}:
\begin{align*}
\dim T_{X}^{1}  &  =(e-4)+\dim T_{\widehat{X}}^{1}+c(X)\\
\dim T_{X}^{2}  &  =(e-2)(e-4)+\dim T_{\widehat{X}}^{2}+c(X).
\end{align*}
The invariant $c(X)$ vanishes for large classes of rational surface
singularities. In \cite{js:rss} de Jong and van Straten proved that $c(X)=0$
for all rational surface singularities with reduced fundamental cycle. This
also follows by the methods \cite{cg:01}. In \cite{ro:92} de Jong proved that
$c(X)=0$ for determinantal rational surface singularities. In the present
paper we generalize this result by using the methods developed in \cite{cg:01}.

\subsection{Some important sheaves and sequences on $\widehat{X}$}

First, let $Y$ $\rightarrow X$ be any modification, and assume that $F\subset
Y$ is a (possible non-reduced) curve which is given by an invertible sheaf of
ideals $\mathcal{I}_{F}$. Let $\mathcal{D}er_{F}(\widehat{X})$ be the subsheaf
of $\Theta_{\widehat{X}}$ consisting of derivations which take $\mathcal{I}%
_{F}$ to itself. Define ${\mathcal{A}}_{F/Y}^{1}$ to be the cokernel of the
map $\Theta_{\widehat{X}}\rightarrow{\mathcal{O}}_{F}(F)$ defined locally --
where $F$ is defined by $x$ -- as $D\mapsto D(x)\otimes\frac{1}{x}\mod (x)$.
Notice that there is an exact sequence
\begin{equation}
0\rightarrow\mathcal{D}er_{F}(F)\rightarrow\Theta_{Y}\rightarrow{\mathcal{O}%
}_{F}(F)\rightarrow{\mathcal{A}}_{F/{Y}}^{1}\rightarrow0 \label{A-exseq}%
\end{equation}
with the maps as above. Similarly, we denote by $\mathcal{T}_{F\subset Y}^{1}$
the cokernel of the map $\mathcal{D}er(\mathcal{O}_{Y},\mathcal{O}%
_{F})\rightarrow{\mathcal{O}}_{F}(F)$ .

In the case $Y=\widehat{X}$ and $F=C$, we have from Proposition 2.3 in
\cite{cg:01} the important exact sequence%
\begin{equation}
0\rightarrow{\mathcal{A}}_{C/{\widehat{X}}}^{1}(C)\rightarrow{\mathfrak{m}%
}{\mathcal{F}}^{1}\rightarrow{\mathcal{T}}_{\widehat{X}}^{1}(C)\rightarrow0.
\label{identities}%
\end{equation}
We remark that the maps in this last sequence are non-canonical. The sequence
sits in a diagram, see \cite[Section 4.2]{cg:01}:
\begin{equation}%
\begin{array}
[c]{ccccccccc}%
0 & \rightarrow & {\mathcal{A}}_{C/{\widehat{X}}}^{1}(C) & \rightarrow &
{\mathfrak{m}}{\mathcal{F}}^{1} & \rightarrow & {\mathcal{T}}_{\widehat{X}%
}^{1}(C) & \rightarrow & 0\\
&  & \downarrow &  & \downarrow &  & \downarrow &  & \\
0 & \rightarrow & \mathcal{T}_{C\subset\widehat{X}}^{1}(C) & \rightarrow &
\mathcal{T}_{C}^{1}(C) & \rightarrow & {\mathcal{T}}_{\widehat{X}}%
^{1}({\mathcal{O}}_{C}(C))(C) & \rightarrow & 0\\
&  & \downarrow &  & \downarrow &  & \downarrow &  & \\
&  & 0 & \rightarrow & {\mathcal{T}}_{\widehat{X}}^{2} & \rightarrow &
{\mathcal{T}}_{\widehat{X}}^{2} & \rightarrow & 0
\end{array}
\label{diagram}%
\end{equation}
The middle horizontal sequence comes from the Zariski-Jacobi long exact
sequence \cite[Th. 5.1]{an:hom} for $\mathbb{C}\rightarrow{\mathcal{O}%
}_{\widehat{X}}(C)\rightarrow{\mathcal{O}}_{C}(C)$ (and is thus natural). The
kernels of the upper vertical maps have support on a finite set of points.

\section{Classification of rational surface singularities with almost reduced
fundamental cycle}

\label{classification}This paper is concerned with a particular class of
rational surface singularities which contains other (more commonly known)
classes of singularities.

\begin{definition}
We say that the fundemental cycle $Z$ is almost reduced if it is reduced on
all non $-2$-curves.
\end{definition}

Rational determinantal surface singularities and more generally rational
quasi-determinantal surface singularities, have almost reduced fundamental
cycle, see \cite[3.4]{wa:equ}, \cite[4.2.1]{ro:92}, \cite{DJ}, \cite{dJ:99}.
Remark also that all two dimensional quotient singularities are
quasi-determinantal and thus have almost reduced fundamental cycle, see
\cite[4.2.2]{ro:92}.

To a normal surface singularity one attach the \textit{dual graph }$\Gamma$ of
the exceptional set in the minimal resolution. The graph $\Gamma$ has a vertex
for each irreducible component $E_{i}$ and there is an edge between two
vertices if the corresponding curves intersect.

In order to understand the shape of the graphs of rational surface
singularities with reduced fundamental cycle, we consider a connected subgraph
$\Gamma_{1}$ of $\Gamma$ which contains only vertices corresponding to
$-2$-curves. We assume further that $\Gamma_{1}$ is maximal in the sense that
all vertices in $\Gamma$ (not in $\Gamma_{1}$) with an edges to a vertex in
$\Gamma_{1}$ corresponds to non $-2$-curves. The graph $\Gamma_{1}$ will
necessarily be the dual graph of a rational double point, and a subgraph
$\Gamma_{1}$ together with all edges in $\Gamma$ which connect to $\Gamma
_{1},$ will be called a \textit{rational double point configuration} or RDP
configuration. This corresponds to the term \textit{erweiterte }%
$(-2)$\textit{-Konfiguration }in \cite[Section 1.5]{ro:92}. Thus a rational
double point configuration is the graph of a rational double point with some
extra edges attached. If there are $n$ such edges we call the rational double
point configuration an $n$-configuration. Thus the $0$-configurations, are the
rational double points. In addition, it may be checked that we only have $1$-,
$2$- and $3$-configurations, see \cite[1.5.1]{ro:92}.

To find the possible rational double point configurations one starts with one
of the wellknown graphs of the rational double points. Then one investigate to
which vertex one are allowed to attach an edge and still have a rational
singularity. In doing this on may consider a sequence $Z_{0}=E,$ $Z_{1},$
$\dots$ $,$ $Z_{k}=Z$ as in section \ref{rational}. This is called a computing
sequence for $Z$. In order to have a rational singularity one must have that
$E_{i}\cong\mathbb{P}^{1}$ for all $i$ and that $Z_{k}\cdot E_{i}>0$ implies
$Z_{k}\cdot E_{i}=1$ for all $k$ in any computation sequence for $Z.$ In
figure \ref{RDPfigure} we give the possible configurations. The configurations
divides into the types $\mathbf{A}$, $\mathbf{D}$ and $\mathbf{E}.$ The
numbers under each vertex is the multiplicity of the fundamental cycle $Z$ at
the corresponding curves. The white vertices correspond to the exceptional
curves $E_{i}$ such that $Z\cdot E_{i}=0.$ The black vertices corresponds to
$E_{i}$ such that $Z\cdot E_{i}>0.$ The subscript of the symbols $\mathbf{A}$,
$\mathbf{D}$ and $\mathbf{E}$ attached to each configuration, give the number
of vertices. For the type $\mathbf{A}$ the superscript $q$ essentially
determines which interior vertex that has an edge attached. Note that for some
values of $(n,q)$ it is allowed for the black vertex to have attached edges.%

\setlength{\unitlength}{1mm}
\newcommand{\unfilledvert}[1]{
\begin{picture}(1,7)\thicklines\circle{2}{\makebox(-3.8,-2)[t]{\tiny$#1$}}
\end{picture}}
\newcommand{\filledvert}[1]{
\begin{picture}(1,7)\thicklines\circle*{2.3}\put(-0.5,0){\makebox
(-3.8,-2)[t]{\tiny$#1$}}
\end{picture}}\newcommand{\solidedge}[0]{
\begin{picture}(4,7)\thicklines\put(-1,0){\line(1,0){5.1}}
\end{picture}}
\newcommand{\dashededge}[0]{
\begin{picture}(4,7)\thicklines\multiput(-1,0)(0.6,0){9}{\line(1,0){0.3}}
\end{picture}}
\newcommand{\upedge}[0]{
\begin{picture}(-1.2,-7)\thicklines\put(-2.1,1){\line(0,1){5}}
\end{picture}}
\newcommand{\upedgevert}[1]{
\begin{picture}(-1.2,-7)\thicklines{\makebox(-4,10.3)[t]{\tiny$#1$}}
\put(-2,7){\circle{2}}
\put(-2.1,1){\line(0,1){5}}
\end{picture}}
\newcommand{\upedgesolidvert}[1]{
\begin{picture}(-1.2,-7)\thicklines{\makebox(-4,10.3)[t]{\tiny$#1$}}
\put(-2,7){\circle*{2.3}}
\put(-2.1,1){\line(0,1){5}}
\end{picture}}
\begin{figure}[h]
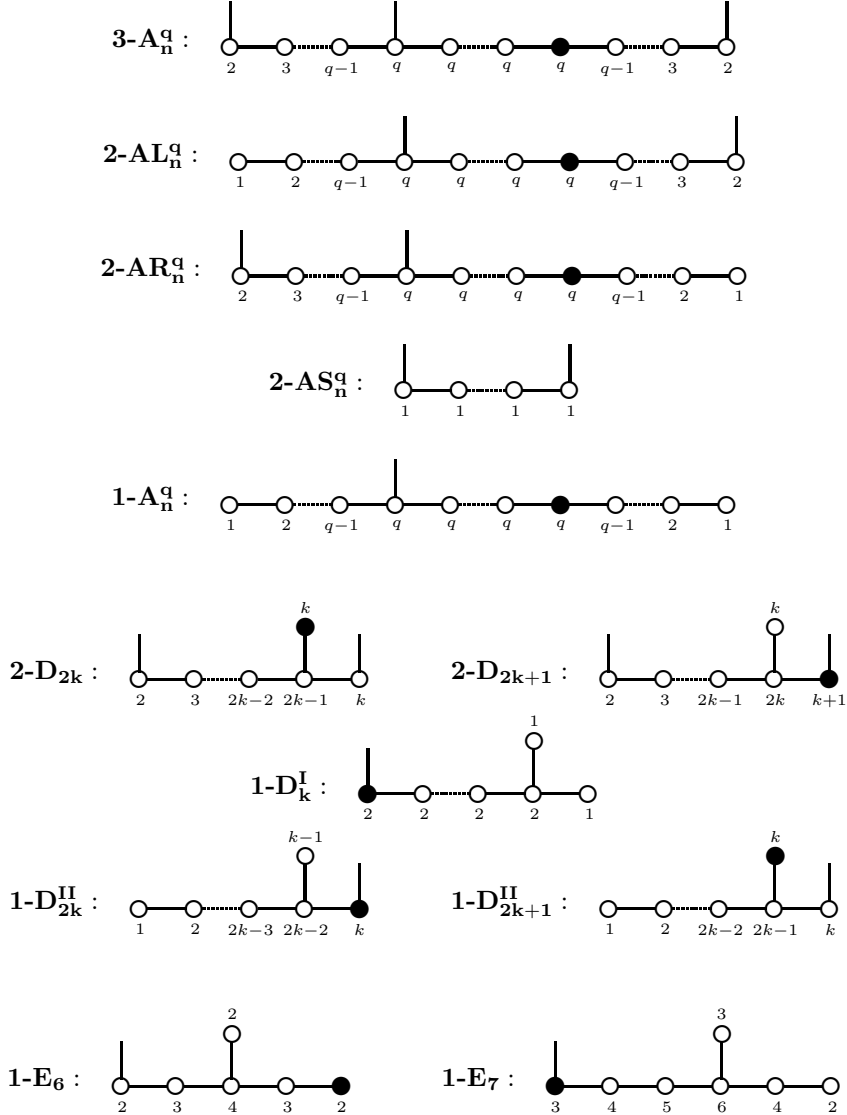

\begin{center}
$\mathbf{3}\mbox{\textbf{-}}\mathbf{A^q_n}:$\hspace{0.3cm}
\unfilledvert{2}\upedge\solidedge\unfilledvert{3}\dashededge\unfilledvert
{q{-}1}\solidedge\unfilledvert{{q}}\upedge\solidedge\unfilledvert
{q}\dashededge\unfilledvert{q}\solidedge\filledvert{q}\solidedge
\unfilledvert{q{-}1}\dashededge\unfilledvert{3}\solidedge\unfilledvert
{2}\upedge\\ \vspace{0.7cm}
$\mathbf{2}\mbox{\textbf{-}}\mathbf{AL^q_n}:$\hspace{0.3cm}
\unfilledvert{1}\solidedge\unfilledvert{2}\dashededge\unfilledvert{q{-}%
1}\solidedge\unfilledvert{q}\upedge\solidedge\unfilledvert{q}\dashededge
\unfilledvert{q}\solidedge\filledvert{q}\solidedge\unfilledvert{q{-}%
1}\dashededge\unfilledvert{3}\solidedge\unfilledvert{2}\upedge\\ \vspace
{0.7cm}
$\mathbf{2}\mbox{\textbf{-}}\mathbf{AR^q_n}:$\hspace{0.3cm}
\unfilledvert{2}\upedge\solidedge\unfilledvert{3}\dashededge\unfilledvert
{q{-}1}\solidedge\unfilledvert{q}\upedge\solidedge\unfilledvert{q}%
\dashededge\unfilledvert{q}\solidedge\filledvert{q}\solidedge\unfilledvert
{q{-}1}\dashededge\unfilledvert{2}\solidedge\unfilledvert{1}\\ \vspace{0.7cm}
$\mathbf{2}\mbox{\textbf{-}}\mathbf{AS^q_n}:$\hspace{0.3cm}
\unfilledvert{1}\upedge\solidedge\unfilledvert{1}\dashededge\unfilledvert
{1}\solidedge\unfilledvert{1}\upedge\\ \vspace{0.7cm}
$\mathbf{1}\mbox{\textbf{-}}\mathbf{A^q_n}:$\hspace{0.3cm}
\unfilledvert{1}\solidedge\unfilledvert{2}\dashededge\unfilledvert{q{-}%
1}\solidedge\unfilledvert{q}\upedge\solidedge\unfilledvert{q}\dashededge
\unfilledvert{q}\solidedge\filledvert{q}\solidedge\unfilledvert{q{-}%
1}\dashededge\unfilledvert{2}\solidedge\unfilledvert{1}
\\ \vspace{1.5cm}
$\mathbf{2}\mbox{\textbf{-}}\mathbf{D_{2k}}:$\hspace{0.3cm}
\unfilledvert{2}\upedge\solidedge\unfilledvert{3}\dashededge\unfilledvert
{2k{-}2}\solidedge\unfilledvert{2k{-}1}\upedgesolidvert{k}\solidedge
\unfilledvert{k}\upedge\hspace{1cm}
$\mathbf{2}\mbox{\textbf{-}}\mathbf{D_{2k+1}}:$\hspace{0.3cm}
\unfilledvert{2}\upedge\solidedge\unfilledvert{3}\dashededge\unfilledvert
{2k{-}1}\solidedge\unfilledvert{2k}\upedgevert{k}\solidedge\filledvert{k{+}%
1}\upedge\\ \vspace{0.7cm}
$\mathbf{1}\mbox{\textbf{-}}\mathbf{D_{k}^I}:$\hspace{0.3cm}
\filledvert{2}\upedge\solidedge\unfilledvert{2}\dashededge\unfilledvert
{2}\solidedge\unfilledvert{2}\upedgevert{1}\solidedge\unfilledvert{1}
\\ \vspace{0.7cm}
$\mathbf{1}\mbox{\textbf{-}}\mathbf{D_{2k}^{II}}:$\hspace{0.3cm}
\unfilledvert{1}\solidedge\unfilledvert{2}\dashededge\unfilledvert{2k{-}%
3}\solidedge\unfilledvert{2k{-}2}\upedgevert{k{-}1}\solidedge\filledvert
{k}\upedge\hspace{1cm}
$\mathbf{1}\mbox{\textbf{-}}\mathbf{D_{2k+1}^{II}}:$\hspace{0.3cm}
\unfilledvert{1}\solidedge\unfilledvert{2}\dashededge\unfilledvert{2k{-2}%
}\solidedge\unfilledvert{2k{-}1}\upedgesolidvert{k}\solidedge\unfilledvert
{k}\upedge\\ \vspace{1.5cm}
$\mathbf{1}\mbox{\textbf{-}}\mathbf{E_6}:$\hspace{0.3cm}
\unfilledvert{2}\upedge\solidedge\unfilledvert{3}\solidedge\unfilledvert
{4}\upedgevert{2}\solidedge\unfilledvert{3}\solidedge\filledvert{2}
\hspace{1cm}
$\mathbf{1}\mbox{\textbf{-}}\mathbf{E_7}:$\hspace{0.3cm}
\filledvert{3}\upedge\solidedge\unfilledvert{4}\solidedge\unfilledvert
{5}\solidedge\unfilledvert{6}\upedgevert{3}\solidedge\unfilledvert
{4}\solidedge\unfilledvert{2}
\end{center}
\caption{\label{RDPfigure}%
Classification of rational doubel point configurations.}
\end{figure}%

We refer to \cite[1.5.1]{ro:92} for the proof that figure \ref{RDPfigure}
gives all possible RDP configurations:

\begin{proposition}
\label{badconfig} The only possible rational double point configurations are
the ones which are given in figure \ref{RDPfigure}.
\end{proposition}

\section{The proof of the main theorem}

We assume that $X$ is a singularity with almost reduced fundamental cycle $Z$.
In order to compute $c(X)=h^{1}(\widehat{X},\mathfrak{m}\mathcal{F}^{1})$ we
use the sequence (\ref{identities}). Taking cohomology we get%
\[
H^{0}(\mathcal{T}_{\widehat{X}}^{1}(C))\overset{\beta}{\longrightarrow}%
H^{1}(\mathcal{A}_{C/\widehat{X}}^{1}(C))\longrightarrow H^{1}(\mathfrak{m}%
\mathcal{F}^{1})\rightarrow0.
\]
In principle one may compute $H^{1}(\mathfrak{m}\mathcal{F}^{1})$ as the
cokernel of $\beta,$ but since this would involve a computation of
$\mathcal{T}_{\widehat{X}}^{1},$ our idea is to consider only the part of
$\mathcal{T}_{\widehat{X}}^{1}$ corresponding to deformations which come from
$\widetilde{X}.$ To do this, we consider the composed map $\eta:H^{1}%
(\widetilde{X},\Theta_{\widetilde{X}})\rightarrow H^{0}(\mathcal{T}%
_{\widehat{X}}^{1})\cong H^{0}(\mathcal{T}_{\widehat{X}}^{1}(C))$ where the
first map is the blow down map composed with the restriction to neigborhoods
of the singular points, see for instance \cite[Th 1.4 and 1.5, 1.6]{wa:equi},
and where the isomorphism comes from the fact that the sheaf $\mathcal{T}%
_{\widehat{X}}^{1}$ has support on points. Composing $\eta$ with $\beta$, we
get a map $\alpha:H^{1}(\widetilde{X},\Theta_{\widetilde{X}})\rightarrow
H^{0}(\widehat{X},\mathcal{A}_{C/\widehat{X}}^{1}(C))$ and a surjection
$\operatorname*{coker}\alpha\twoheadrightarrow H^{1}(\mathfrak{m}%
\mathcal{F}^{1}).$ We then proceed by showing that $\alpha$ may be computed
separately for each RDP-configuration.

The main parts of the proof is divided into to lemmas. In the first lemma, we
compute $H^{1}(\mathcal{A}_{C/\widehat{X}}^{1}(C))$, (for technical reasons in
the proof of the second lemma we need a slightly more general statement), and
in the second lemma we proceed to compute $\operatorname*{coker}\alpha.$

\begin{lemma}
\label{lemme:iii}Let $I$ be such that $i\in I$ implies $Z\cdot E_{i}=0.$ Let
$Y$ be obtained from $\widetilde{X}$ by contracting $E_{i}$ for $i\in I.$ Then
$\widetilde{\pi}:\widetilde{X}\rightarrow\widehat{X}$ factors through
$\widetilde{\pi}_{Y}:\widetilde{X}\rightarrow Y$ and $F=(\widetilde{\pi}%
_{Y})_{\ast}Z$ is given by an invertible sheaf of ideals on $Y.$ We have that
$H^{1}(\mathcal{A}_{F/Y}^{1}(F))\cong H^{1}(\widetilde{X},\mathcal{O}%
_{Z-E}(2Z))\cong H^{1}(\widehat{X},\widetilde{\pi}_{\ast}\mathcal{O}%
_{Z-E}(2Z)).$ Moreover, if $K$ is the canonical divisor on $\widetilde{X}$
\ we have $h^{1}(\mathcal{A}_{F/Y}^{1}(F))=(Z-E)\cdot(K-Z)=\sum(s_{i}-1)$
where the sum is taken over all RDP-configurations and where $s_{i}$ is the
weight of the black vertex (or $s_{i}=1$ if there is not a black vertex) in
the corresponding configuration (see figure \ref{RDPfigure}).
\end{lemma}

\begin{proof}
The factorisation of $\widetilde{\pi}$ follows from Theorem \ref{Tyu}, and the
ideal sheaf of $F$ is $\mathfrak{m}\mathcal{O}_{Y}$ and is invertible.
Consider the diagram of exact sequences%
\begin{equation}%
\begin{array}
[c]{ccccccc}%
0 &  & 0 &  &  &  & \\
\downarrow &  & \downarrow &  &  &  & \\
H^{1}(Y,\Theta_{Y}(F)) & \rightarrow & H^{1}(\mathcal{O}_{F}(2F)) &
\rightarrow & H^{1}(\mathcal{A}_{F/\widehat{X}}^{1}(F)) & \rightarrow & 0\\
\downarrow &  & \text{ }\downarrow\cong &  & \text{ \ }\downarrow\gamma &  &
\\
H^{1}(\Theta_{\widetilde{X}}(Z)) & \rightarrow & H^{1}(\mathcal{O}_{Z}(2Z)) &
\rightarrow & H^{1}(\mathcal{T}_{Z}^{1}(Z)) & \rightarrow & 0\\
\downarrow &  & \downarrow &  & \downarrow &  & \\
H^{0}(R^{1}(\widetilde{\pi}_{Y})_{\ast}\Theta_{\widetilde{X}}(Z)) &
\rightarrow & 0 &  & 0 &  &
\end{array}
\label{thediagram}%
\end{equation}
The first row comes from (\ref{A-exseq}) and the second row comes from the
sequence%
\begin{equation}
0\rightarrow\mathcal{D}er_{Z}(\mathcal{O}_{\widetilde{X}}(Z))\rightarrow
\Theta_{\widetilde{X}}(Z)\rightarrow\mathcal{O}_{Z}(2Z)\rightarrow
\mathcal{T}_{Z}^{1}(Z)\rightarrow0. \label{four-term}%
\end{equation}
The two left-most vertical rows come from the Lerray spectral sequence using
that $(\widetilde{\pi}_{Y})_{\ast}\Theta_{\widetilde{X}}=\Theta_{Y},$ (see
\cite{bw:loc}) $(\widetilde{\pi}_{Y})_{\ast}\mathcal{O}_{Z}=\mathcal{O}_{F}$
and $R^{1}(\widetilde{\pi}_{Y})_{\ast}\mathcal{O}_{Z}(nZ)=0,$ which follows
since only curves $E_{i}$ such that $Z\cdot E_{i}=0$ are contracted in $Y.$ A
local calculation shows that $\mathcal{D}er_{Z}(\mathcal{O}_{\widetilde{X}%
})=\Theta_{\widetilde{X}}(\log E)=:S,$ and as in \cite[2.2]{wa:equi} we get
that (\ref{four-term}) splits in two short exact sequences
\begin{equation}
0\rightarrow S\rightarrow\Theta_{\widetilde{X}}(Z)\rightarrow\oplus
\mathcal{O}_{E_{i}}(E_{i}+Z)\rightarrow0 \label{first-sequence}%
\end{equation}
and%
\begin{equation}
0\rightarrow\oplus\mathcal{O}_{E_{i}}(E_{i}+Z)\rightarrow\mathcal{O}%
_{Z}(2Z)\rightarrow\mathcal{T}_{Z}^{1}(Z)\rightarrow0. \label{last-seq}%
\end{equation}
From the long exact sequences in cohomology of this last sequence and from
diagram (\ref{thediagram}) (using the Snake lemma) we get that
\[
\ker\gamma\cong\frac{H^{1}(\oplus\mathcal{O}_{E_{i}}(E_{i}+Z))}%
{\operatorname*{im}H^{0}(\mathcal{T}_{Z}^{1}(Z))+\operatorname*{im}%
H^{1}(Y,\Theta_{Y}(F))}.
\]
Using the cohomology sequence of (\ref{first-sequence}) from the diagram
(\ref{thediagram}) we get%
\[
\frac{H^{1}(\oplus\mathcal{O}_{E_{i}}(E_{i}+Z))}{\operatorname*{im}%
H^{1}(Y,\Theta_{Y}(F))}\cong\frac{H^{1}(\Theta_{\widetilde{X}}(Z))}%
{\operatorname*{im}H^{1}(S)+\operatorname*{im}H^{1}(Y,\Theta_{Y}(F))}%
\cong\frac{H^{0}(R^{1}(\widetilde{\pi}_{Y})_{\ast}\Theta_{\widetilde{X}}%
(Z))}{\operatorname*{im}H^{0}(R^{1}(\widetilde{\pi}_{Y})_{\ast}S)}%
\]
The sequence (\ref{first-sequence}) gives
\[
0\rightarrow R^{1}(\widetilde{\pi}_{Y})_{\ast}S\rightarrow R^{1}%
(\widetilde{\pi}_{Y})_{\ast}\Theta_{\widetilde{X}}(Z)\rightarrow
R^{1}(\widetilde{\pi}_{Y})_{\ast}(\oplus\mathcal{O}_{E_{i}}(E_{i}%
+Z))\rightarrow0.
\]
These sheaves are supported on the singular points of $\widehat{X},$ so we
get
\[
\frac{H^{0}(R^{1}(\widetilde{\pi}_{Y})_{\ast}\Theta_{\widetilde{X}}%
(Z))}{\operatorname*{im}H^{0}(R^{1}(\widetilde{\pi}_{Y})_{\ast}S)}\cong
H^{0}(R^{1}(\widetilde{\pi}_{Y})_{\ast}(\oplus\mathcal{O}_{E_{i}}%
(E_{i}+Z)))\cong\oplus_{i\in I}H^{1}(\mathcal{O}_{E_{i}}(E_{i}))
\]
and%
\[
\ker\gamma\cong\frac{\oplus_{i\in I}H^{1}(\mathcal{O}_{E_{i}}(E_{i}%
))}{\operatorname*{im}H^{0}(\mathcal{T}_{Z}^{1}(Z))}.
\]
To find $\ker\gamma$ we thus need to find the image of $H^{0}(\mathcal{T}%
_{Z}^{1}(Z))$ in $\oplus_{i\in I}H^{1}(\mathcal{O}_{E_{i}}(E_{i})).$ There is
an (unnatural) exact sequence
\begin{equation}
0\rightarrow\mathcal{T}_{E}^{1}\rightarrow\mathcal{T}_{Z}^{1}(Z)\rightarrow
\mathcal{O}_{Z-E}(2Z)\rightarrow0, \label{wahl-sec}%
\end{equation}
see \cite[(2.6.3)]{wa:sim}. As in the proof of proposition 4.3 in
\cite{cg:01}, we have $H^{0}(\mathcal{O}_{Z-E}(2Z))=0.$ Thus $H^{0}%
(\mathcal{T}_{Z}^{1}(Z))=H^{0}(\mathcal{T}_{E}^{1})=\mathbb{C}^{s}$ where $s$
is the number of intersection points in $E.$ The map $H^{0}(\mathcal{T}%
_{Z}^{1}(Z))\rightarrow\oplus_{i\in I}H^{1}(\mathcal{O}_{E_{i}}(E_{i}))$ is
the composition of the (only) connecting morphism $H^{0}(\mathcal{T}_{Z}%
^{1}(Z))\rightarrow H^{1}(\oplus\mathcal{O}_{E_{i}}(E_{i}+Z))$ (resulting from
\ref{last-seq}) with the projection $H^{1}(\oplus\mathcal{O}_{E_{i}}%
(E_{i}+Z))\rightarrow\oplus_{i\in I}H^{1}(\mathcal{O}_{E_{i}}(E_{i})).$ Using
(\ref{wahl-sec}) we may compute this by the connecting morphism from
$0\rightarrow\oplus\mathcal{O}_{E_{i}}(E_{i})\rightarrow\mathcal{O}%
_{E}(E)\rightarrow\mathcal{T}_{E}^{1}\rightarrow0$ composed with the
projection. If we consider $r,$ $-2$-curves $E_{i}$ intersecting a $-m$-curve
$E_{j}$, we claim that the map $H^{0}(\mathcal{T}_{E_{j}\cup(\cup E_{i})}%
^{1})=\mathbb{C}^{r}\rightarrow\oplus H^{1}(\mathcal{O}_{E_{i}}(E_{i}))\oplus
H^{1}(\mathcal{O}_{Ej}(E_{j}))=\mathbb{C}^{r}\oplus\mathbb{C}^{m-1}$ is given
by $e_{i}\mapsto(e_{i},\eta_{i})$ where the $\eta_{i}$ are linearly
independent if $r<m$ and span $\mathbb{C}^{m-1}$ if $r\geq m-1.$ This can be
seen by a tedious but straight forward calculation. Using this we can deduce
that $H^{0}(\mathcal{T}_{E}^{1})\rightarrow\oplus_{i\in I}H^{1}(\mathcal{O}%
_{E_{i}}(E_{i}))$ is surjective (and hence $\ker\gamma=0$) in the present
situation by checking the different RDP-configurations in figure
\ref{RDPfigure}: For each connected set of white vertices, there is a vertex
(corresponding to $E_{i}$) with an edge to a black vertex (corresponding to
$E_{j}.$) Since $H^{0}(\mathcal{T}_{E_{j}\cup E_{i}}^{1})=\mathbb{C}%
\rightarrow H^{1}(\mathcal{O}_{E_{i}}(E_{i}))\oplus H^{1}(\mathcal{O}%
_{Ej}(E_{j}))=\mathbb{C}\oplus\mathbb{C\rightarrow}H^{1}(\mathcal{O}_{E_{i}%
}(E_{i}))=\mathbb{C}$ maps the generator to a generator, we will have the
corresponding basis vector $e_{i}$ $\in\operatorname*{im}H^{0}(\mathcal{T}%
_{E}^{1})\subset\oplus_{i\in I}H^{1}(\mathcal{O}_{E_{i}}(E_{i})).$ If we have
two $-2$-curves $E_{l}$ and $E_{k}$ such that $Z\cdot E_{l}=Z\cdot E_{k}=0$
that intersect, we will likewise have $e_{l}+e_{k}\in\operatorname*{im}%
H^{0}(\mathcal{T}_{E}^{1}).$ From this follows that we have $e_{i}%
\in\operatorname*{im}H^{0}(\mathcal{T}_{E}^{1})$ for all $-2$-curves $E_{i}$
in $\oplus_{i\in I}H^{1}(\mathcal{O}_{E_{i}}(E_{i})).$ If there is a
$-m$-curve $E_{j}$ (where $Z$ is reduced) which intersects $Z$ in zero there
must be $m$ curves which intersect this curve. From this we now see the
corresponding $\eta_{i}$ spanning $H^{1}(\mathcal{O}_{Ej}(E_{j}))$ must be in
$\operatorname*{im}H^{0}(\mathcal{T}_{E}^{1}),$ and thus that $H^{0}%
(\mathcal{T}_{E}^{1})\rightarrow\oplus_{i\in I}H^{1}(\mathcal{O}_{E_{i}}%
(E_{i}))$ must be surjective. We get that $\ker\gamma=0$ and it follows that
$\gamma$ is an isomorphism. \ 

Taking cohomology of the sequence (\ref{wahl-sec}), we get $H^{1}%
(\mathcal{T}_{Z}^{1}(Z))\cong H^{1}(\widetilde{X},\mathcal{O}_{Z-E}(2Z)).$
Since $R^{1}(\widetilde{\pi}_{Y})_{\ast}\mathcal{O}_{Z-E}(2Z)=0$ we also get
$H^{1}(\widetilde{X},\mathcal{O}_{Z-E}(2Z))\cong H^{1}(\widehat{X}%
,(\widetilde{\pi}_{Y})_{\ast}\mathcal{O}_{Z-E}(2Z)).$

Finally, we have $h^{1}(\mathcal{T}_{Z}^{1}(Z))=h^{1}(\mathcal{O}%
_{Z-E}(2Z))=(Z-E)\cdot(K-Z)$ where the last formula follows as in the proof of
proposition 4.3 in \cite{cg:01}. This allow us to compute $s-1$ as the
contribution from each RDP-configuration.
\end{proof}

By the lemma, $H^{1}(\mathcal{A}_{C/\widehat{X}}^{1}(C))\cong H^{1}%
(\widetilde{X},\mathcal{O}_{Z-E}(2Z)),$ and we may view $\alpha$ as a map
$H^{1}(\widetilde{X},\Theta)\rightarrow H^{1}(\widetilde{X},\mathcal{O}%
_{Z-E}(2Z)).$ The main technical observation regarding this map is given in
the following lemma.

\begin{lemma}
\label{alphalemma}Let $\widetilde{U}_{RDP,i}$ be (small enough) neighborhoods
of (the $-2$-curves of) each RDP-configuration in $\widetilde{X}.$ The map
$\alpha:H^{1}(\widetilde{X},\Theta)\rightarrow H^{1}(\widetilde{X}%
,\mathcal{O}_{Z-E}(2Z))$ is the composition of the restriction map
$H^{1}(\widetilde{X},\Theta_{\widetilde{X}})\rightarrow\oplus H^{1}%
(\widetilde{U}_{RDP,i},\Theta_{\widetilde{X}})$ with a direct sum of maps
$\alpha_{i}:H^{1}(\widetilde{U}_{RDP,i},\Theta_{\widetilde{X}})\rightarrow
H^{1}(\widetilde{U}_{RDP,i},\mathcal{O}_{Z-E}(2Z)).$ Moreover, these last maps
are surjective except in the $\mathbf{3}$\textbf{-}$\mathbf{A}_{\mathbf{n}}%
$-case where the cokernel is $\mathbb{C.}$
\end{lemma}

\begin{proof}
The proof is divided in two steps.

\textbf{Step 1: }We prove that we may reduce to the case where $\widehat{X}$
has only rational double points. Define $Y$ to be the surface obtained from
$\widetilde{X}$ by blowing down all $-2$-curves $E_{i}$ such that $Z\cdot
E_{i}=0.$ Then $\widetilde{\pi}$ factors $\widetilde{X}\overset{\widetilde
{\pi}_{Y}}{\rightarrow}Y\overset{\pi_{Y}}{\rightarrow}\widehat{X}.$ We put
$F=(\widetilde{\pi}_{Y})_{\ast}Z.$ From the Zariski-Jacobi long exact sequence
for $\mathbb{C}\rightarrow{\mathcal{O}}_{Y}\rightarrow{\mathcal{O}}_{F}$ we
get a sequence%
\[
0\rightarrow\mathcal{T}_{F\subset Y}^{1}(F)\rightarrow\mathcal{T}_{F}%
^{1}(F)\rightarrow\mathcal{T}_{Y}^{1}(\mathcal{O}_{F}(F))\rightarrow0
\]
and a composed map%
\[
H^{1}(\widetilde{X},\Theta)\rightarrow H^{0}(\mathcal{T}_{Y}^{1})\cong
H^{0}(\mathcal{T}_{Y}^{1}(F))\rightarrow H^{1}(\mathcal{T}_{F\subset Y}%
^{1}(F))
\]
There is a surjection $\mathcal{A}_{F/Y}^{1}\rightarrow\mathcal{T}_{F\subset
Y}^{1}$ and the kernel is supported on a finite set of points. From this we
get $H^{1}(\mathcal{A}_{F/Y}^{1}(F))\cong H^{1}(\mathcal{T}_{F\subset Y}%
^{1}(F))$ and from lemma \ref{lemme:iii} this is (canonically) isomorphic to
$H^{1}(\mathcal{O}_{Z-E}(2Z)).$ Thus we get a map $\alpha^{\prime}%
:H^{1}(\widetilde{X},\Theta)\rightarrow H^{1}(\mathcal{O}_{Z-E}(2Z)).$ We
claim that $\alpha^{\prime}=\alpha.$ To see this, we compute the map in
\v{C}ech-cohomology. Cover $\widehat{X}$ by open affines $U_{i}$ such that
$U_{i}\cap U_{j}$ are small enough and in particular do not contain
singularities. Note that the part of $F$ which is contracted in $\widehat{X}$
is reduced (away from the singularities.) From this follows that $H^{1}%
(U_{i}\times_{\widehat{X}}Y,\mathcal{T}_{F}^{1}(F))=H^{1}(U_{i}\times
_{\widehat{X}}Y,\mathcal{T}_{F\subset Y}^{1}(F))=0,$ and we may use $U_{i}%
^{Y}:=U_{i}\times_{\widehat{X}}Y$ as covering of $Y.$ Let $\xi\in
H^{1}(\widetilde{X},\Theta_{\widetilde{X}})$ . Denoting by $\mathcal{N}_{F/Y}$
the normal sheaf of $F$ in $Y$, we have $H^{1}(U_{i}^{Y},\mathcal{N}%
_{F/Y})=H^{1}(U_{i}^{Y},\mathcal{O}_{F}(F))=0$ since $F$ must be principal on
$U_{i}^{Y}.$ Setting $F_{i}=F\cap U_{i}^{Y},$ we also have $H^{0}(U_{i}%
^{Y},\mathcal{T}_{F/Y}^{2})=0$ and from the exact sequence $0\rightarrow
H^{1}(U_{i}^{Y},\mathcal{N}_{F/Y})\rightarrow T_{F_{i}/U_{i}^{Y}}%
^{2}\rightarrow H^{0}(U_{i}^{Y},\mathcal{T}_{F/Y}^{2})\rightarrow0$ we have
$T_{F_{i}/U_{i}^{Y}}^{2}=0.$ From this follows that for each $i$ there are
deformations $\overline{F}_{i}$ (over $\mathbb{C[\varepsilon]}$) of $F_{i}$
and inclusions $\overline{F}_{i}\subset\overline{U}_{i}^{Y}$ where
$\overline{U}_{i}^{Y}$ represents the image of $\xi$ under the blow down and
restriction map. We may blow down $\overline{F}_{i}$ to a deformation
$\overline{C}_{i}$ of $C_{i}=C\cap U_{i}.$ Denote by $\nu_{i}$ the element
corresponding to $\overline{F}_{i}$ in $H^{0}(U_{i}^{Y},$ $\mathcal{T}_{F}%
^{1}(F))\cong H^{0}(U_{i}^{Y},\mathcal{T}_{F}^{1})$ ($F$ is principal on
$U_{i}^{Y}$), and denote by $\eta_{i}$ the element corresponding to
$\overline{C}_{i}$ in $H^{0}(U_{i},\mathcal{T}_{C}^{1}(C))\cong H^{0}(U_{i},$
$\mathcal{T}_{C}^{1}).$ Then $\nu_{i}$ and $\xi$ will have the same image in
$H^{0}(U_{i}^{Y},$ $\mathcal{T}_{Y}^{1}(\mathcal{O}_{F}(D)))$ and $\eta_{i}$
and $\xi$ will have the same image in $H^{0}(U_{i},$ $\mathcal{T}_{\widehat
{X}}^{1}(\mathcal{O}_{C}(C))).$ Moreover the restriction of $\nu_{i}$ and
$\eta_{i}$ to $U_{i}\cap U_{j}$ $\cong$ $U_{i}^{Y}\cap U_{j}^{Y}$ give the
same element in $H^{0}(U_{i}\cap U_{j},\mathcal{T}_{C}^{1}(C))\cong$
$H^{0}(U_{i}^{Y}\cap U_{j}^{Y},\mathcal{T}_{F}^{1}(F)).$ In particular,
$(\eta_{i})_{|U_{i}\cap U_{j}}-(\eta_{j})_{|U_{i}\cap U_{j}}$ and $(\nu
_{i})_{|U_{i}\cap U_{j}}-(\nu_{j})_{|U_{i}\cap U_{j}}$ give the same element
in $H^{0}(U_{i}\cap U_{j},\mathcal{T}_{C\subset\widehat{X}}^{1}(C)).$ This
shows that $\alpha=\alpha^{\prime}.$ Also, since the map $H^{1}(\widetilde
{X},\Theta)\rightarrow H^{0}(\mathcal{T}_{Y}^{1})$ factors through $\oplus
H^{1}(\widetilde{U}_{RDP,i},\Theta_{\widetilde{X}}),$ so does $\alpha^{\prime
}.$ From the calculation of $\alpha^{\prime}$ above (and the proof of lemma
\ref{lemme:iii}) it further follows that $\oplus H^{1}(\widetilde{U}%
_{RDP,i},\Theta_{\widetilde{X}})\rightarrow H^{1}(\widetilde{X},\mathcal{O}%
_{Z-E}(2Z))$ splits into a sum as stated. Consider a $-m_{i}$-curve $E_{i}$ in
$\widetilde{X}$ with $m_{i}\neq2.$ If $E_{i}\cdot Z=0,$ we 'change' $Y$ (and
$\widetilde{X}$) in order to increase $m_{i}$ by one. This may be done by
'plumbing' or by considering the blow up in a point on $E_{i}.$ From above it
is clear that this does not change the image of $\alpha.$ Thus we may assume
$E_{i}\cdot Z<0,$ and hence that $Y=\widehat{X}$ have only rational double points.

\textbf{Step 2:} We calculate the cokernel of $\alpha_{i}$ assuming
$E_{i}\cdot Z<0$ for all non-$-2$-curves $E_{i}.$ Let $\widehat{U}_{RDP,i}$
denote the image of $\widetilde{U}_{RDP,i}$ in $\widehat{X}.$ Since
$\widehat{X}$ contains only rational double points, $H^{1}(\widetilde
{U}_{RDP,i},\Theta_{\widetilde{X}})$ surjects to $H^{0}(\widehat{U}%
_{RDP,i},\mathcal{T}_{\widehat{X}}^{1}).$ From this it follows that
$\operatorname*{coker}\alpha=H^{1}(\widehat{X},\mathfrak{m}\mathcal{F}%
^{1})\cong H^{1}(\widehat{X},\mathcal{T}_{C}^{1}(C)).$ The last isomorphism
comes from diagram (\ref{diagram}). We have that $H^{1}(\widehat
{X},\mathcal{T}_{C}^{1}(C))\cong\oplus H^{1}(\widehat{U}_{RDP,i}%
,\mathcal{T}_{C}^{1}(C))$ because the support of $\mathcal{T}_{C}^{1}(C)$ is
contained in $\cup\widehat{U}_{RDP,i}$ since $\mathcal{T}_{C}^{1}(C)$ is zero
at a non-singular points where $C$ is reduced. We thus obtain
$\operatorname*{coker}\alpha_{i}$ as $H^{1}(\widehat{U}_{RDP,i},\mathcal{T}%
_{C}^{1}(C)).$ The calculation divides in cases for the different rational
double point configurations. For the $\mathbf{2}$-$\mathbf{AS}_{\mathbf{n}%
}^{\mathbf{q}}$-case the sheaf $\mathcal{T}_{C}^{1}(C)$ restricted to
$\widehat{U}_{RDP,i}$ has support on a finite set of points, so $H^{1}%
(\widehat{U}_{RDP,i},\mathcal{T}_{C}^{1}(C))=0.$ In the other cases there is a
unique $-2$-curve which we denote by $H$, which is not contracted in
$\widehat{X}.$ The image of $H$ in $\widehat{X}$ will contain at least one
singular point $p$ for $\widehat{X}.$ Since our calculation only will depend
on a formal neighborhood of the exceptional set, we may assume that
$U_{i}=\operatorname*{Spec}\mathbb{C}[u_{i},v_{i}]$ cover the exceptional
curve $H.$ We further assume that $\widehat{U}_{RDP,i}=V_{0}\cup V_{1}$ with
$V_{i}=\operatorname*{Spec}B_{i}$ and that $\widetilde{\pi}$ induce
$\varphi_{i}:B_{i}\rightarrow\mathbb{C}[u_{i},v_{i}].$ We also assume that
$U_{0}\cap U_{1}\cong V_{0}\cap V_{1}.$ We calculate $H^{1}(\widehat
{U}_{RDP,i},\mathcal{T}_{C}^{1}(C))$ as
\[
H^{1}(\widehat{U}_{RDP,i},\mathcal{T}_{C}^{1}(C))\cong H^{0}(V_{0}\cap
V_{1},\mathcal{T}_{C}^{1}(C))/H^{0}(V_{0},\mathcal{T}_{C}^{1}(C))+H^{0}%
(V_{1},\mathcal{T}_{C}^{1}(C)).
\]
To find $\mathcal{T}_{C}^{1}(C)$ locally, we consider $B=B_{0}=P/I$ where $P$
is a regular algebra with parameters $z_{1},z_{2},z_{3}$ and $I=(g).$ Assuming
that $x$ defines $C$ locally, we have%
\[
\operatorname*{Der}(P,B)\rightarrow\operatorname*{Hom}\nolimits_{P}%
((g,x),B/(x))\rightarrow H^{0}(V_{0},\mathcal{T}_{C}^{1}(C))\rightarrow0
\]
and similarly on $V_{1}.$ Tere is an isomorphism $H^{0}(U_{0}\cap
U_{1},\mathcal{O}_{Z-E}(2Z))\cong H^{0}(V_{0}\cap V_{1},\mathcal{T}_{C}%
^{1}(C)).$ If $P^{\prime}$ is a localisation of $P$ such that $B^{\prime
}=H^{0}(V_{0}\cap V_{1},\mathcal{O}_{C})\cong H^{0}(U_{0}\cap U_{1}%
,\mathcal{O}_{Z})=P^{\prime}/(g,x),$ the isomorphism is given on the level of
representatives by mapping $h\in B^{\prime}$ to the element $\phi_{h}%
\in\operatorname*{Hom}\nolimits_{P^{\prime}}((g,x),B^{\prime})$ such that
$\phi_{h}(g)=0$ and $\phi_{h}(x)=h.$ This gives a composed map $B^{\prime
}\twoheadrightarrow H^{0}(V_{0}\cap V_{1},\mathcal{T}_{C}^{1}%
(C))\twoheadrightarrow H^{1}(\widehat{U}_{RDP,i},\mathcal{T}_{C}^{1}(C))$
which factors through $H^{1}(\widetilde{U}_{RDP,i},\mathcal{O}_{Z-E}(2Z))$. We
compute the image of this map:

$\mathbf{A}_{\mathbf{n}}$\textbf{-case: }The singularity $p$ will be of type
$A_{r}$ ($r=q+k-1,$ $q-1$ or $q-2$) and we may assume that $\varphi_{i}$ is
given by $z_{1}=u_{i}^{r}v_{i}^{r+1},z_{2}=u_{i},z_{3}=u_{i}v_{i}.$ The curve
$H$ is given by $v_{0}=0$ and $v_{1}=0,$ ($u_{0}=1/u_{1},$ $v_{0}=u_{1}%
^{2}v_{1}$) and the exceptional set for $p=p_{i}$ is given by $u_{i}=0.$ We
have $g=z_{1}z_{2}-z_{3}^{r+1}$. One finds that $u_{0}^{s-1}v_{0}^{s},0\leq
s,$ generates
\[
H^{1}(\widetilde{U}_{RDP,i},\mathcal{O}_{Z-E}(2Z))\cong\frac{\mathbb{C}%
[u_{0},v_{0},u_{0}^{-1}]/(v_{0}^{q-1})}{\operatorname*{im}\mathbb{C}%
[u_{0},u_{0}v_{0},u_{0}^{r}v_{0}^{r+1}]+\operatorname*{im}\frac{1}{u_{0}^{2}%
}\mathbb{C}[u_{0}^{-1},u_{0}v_{0},u_{0}^{t+2}v_{0}^{t+1}]}%
\]
where the singularity $p=p_{0}$ is of type $A_{r}$ and the other singularity
$p_{1}$ in the image of $H$ is of type $A_{t}.$

We have $h:=u_{0}^{s-1}v_{0}^{s}=z_{3}^{s}/z_{2}.$ For the $\mathbf{3}%
$\textbf{-}$\mathbf{A}_{\mathbf{n}}$-case we have $x=z_{1}z_{3},$ we have
$x=z_{1}$ in the $\mathbf{2}$\textbf{-}$\mathbf{AR}_{\mathbf{n}}^{\mathbf{q}}$
case and $x=z_{1}+z_{3}^{q}$ in $\mathbf{2}$\textbf{-}$\mathbf{AL}%
_{\mathbf{n}}^{\mathbf{q}}$ case. In the two latter cases we can define
$D\in\operatorname*{Der}(P,B^{\prime})$ (since we must have $1/z_{2}\in
B^{\prime}$) by $D(z_{1})=$ $z_{3}^{s}/z_{2}$, $D(z_{2})=D(z_{3})=0,$ so
$D(x)=z_{3}^{s}/z_{2}$ and $D(g)=D(z_{1})z_{2}=z_{3}^{s}.$ Define $\psi_{0}%
\in\operatorname*{Hom}\nolimits_{P}((g,x),B/(x))$ by $\psi_{0}(g)=-z_{3}^{s}$
and $\psi_{0}(x)=0.$ Restricted to $V_{0}\cap V_{1}$ we have that the class of
$\psi_{0}$ equals the class of $\phi_{h}.$ This shows that $\phi_{h}$ is in
the image of $H^{0}(V_{0},\mathcal{T}_{C}^{1}(C))$ and that
$\operatorname*{coker}\alpha_{i}=0.$ Similarly, for the $\mathbf{3}$%
\textbf{-}$\mathbf{A}_{\mathbf{n}}$-case we define $D\in\operatorname*{Der}%
(P,B^{\prime})$ by $D(z_{1})=$ $z_{3}^{s-1}/z_{2}$, $D(z_{2})=D(z_{3})=0$
($s\geq1$). This gives $D(x)=z_{3}D(z_{1})=z_{3}^{s}/z_{2}$ and $D(g)=z_{3}%
^{s-1}.$ This shows that $u_{0}^{s-1}v_{0}^{s}\mapsto0$ in $H^{1}(\widehat
{U}_{RDP,i},\mathcal{T}_{C}^{1}(C))$ for $s\geq1.$ On the other hand if
$[\psi_{0}]\in H^{0}(V_{0}\cap V_{1},\mathcal{T}_{C}^{1}(C))\cong
\mathbb{C}[u_{0},v_{0},u_{0}^{-1}]/(v_{0}^{q-1})$ is in the image of $\oplus
H^{0}(V_{i},\mathcal{T}_{C}^{1}(C))$ one finds that it cannot contain the term
$1/u_{0}.$ Thus $1/u_{0}$ cannot map to zero, and we have that
$\operatorname*{coker}\alpha_{i}=\mathbb{C.}$

$\mathbf{D}_{\mathbf{n}}$\textbf{-case:} We will consider the cases
$\mathbf{2}$\textbf{-}$\mathbf{D}_{\mathbf{2k}}$ and $\mathbf{1}$%
\textbf{-}$\mathbf{D}_{\mathbf{2k+1}}^{\mathbf{II}}.$ The cases $\mathbf{2}%
$\textbf{-}$\mathbf{D}_{\mathbf{2k+1}}$ and $\mathbf{1}$\textbf{-}%
$\mathbf{D}_{\mathbf{2k}}^{\mathbf{II}}$ are similar and will be left to the
reader. For the case $\mathbf{1}$\textbf{-}$\mathbf{D}_{\mathbf{k}%
}^{\mathbf{I}}$, see next paragraph. The divisor $H$ is given by $v_{0}=-1$
and the image of $H$ in $\widehat{X}$ contains only one singular point $p$
which will be of type $A_{2k-1}$ in the case $\mathbf{2}$\textbf{-}%
$\mathbf{D}_{\mathbf{2k}}$ and of type $A_{2k}$ in the case $\mathbf{1}%
$\textbf{-}$\mathbf{D}_{\mathbf{2k+1}}^{\mathbf{II}}.$ We may assume that
$\varphi_{0}$ is given by $z_{1}=u_{0}^{2}v_{0},$ $z_{2}=u_{0}^{2k-2}%
v_{0}^{2k-1}$ in the $\mathbf{2}$\textbf{-}$\mathbf{D}_{\mathbf{2k}}$ case and
$z_{2}=u_{0}^{2k-1}v_{0}^{2k}$ in the $\mathbf{1}$\textbf{-}$\mathbf{D}%
_{\mathbf{2k+1}}^{\mathbf{II}},$ and $z_{3}=u_{0}v_{0}.$ The singularity $p$
is given by $g=z_{1}z_{2}-z_{3}^{2k}$ and $g=z_{1}z_{2}-z_{3}^{2k+1}$
respectively. The (reduced) inverse image of $p$ in $U_{0}$ is given by
$u_{0}v_{0}.$ We may assume that $x$ (which defines $C$) is given by%
\[
x=\left\{
\begin{array}
[c]{lllll}%
u_{0}^{2k-1}v_{0}^{k}(v_{0}+1)^{k} & = & z_{2}z_{3}+\sum_{j=0}^{k-1}\binom
{k}{j}z_{1}^{k-j-1}z_{3}^{2j+1} & \text{for} & \mathbf{2}\text{\textbf{-}%
}\mathbf{D}_{\mathbf{2k}}\\
u_{0}^{2k-1}v_{0}^{k}(v_{0}+1)^{k} & = & z_{2}+\sum_{j=0}^{k-1}\binom{k}%
{j}z_{1}^{k-j-1}z_{3}^{2j+1} & \text{for} & \mathbf{1}\text{\textbf{-}%
}\mathbf{D}_{\mathbf{2k+1}}^{\mathbf{II}}%
\end{array}
\right.
\]
Let
\[
H_{t}:=\frac{\mathbb{C}[u_{0},v_{0},u_{0}^{-1}]/((v_{0}+1)^{t})}%
{\operatorname*{im}\mathbb{C}[u_{0}^{2}v_{0},u_{0}v_{0},u_{0}^{r}v_{0}%
^{r+1}]+\operatorname*{im}\frac{1}{u_{0}^{2}}\mathbb{C}[u_{0}^{-1},u_{0}%
^{2}(v_{0}+1)]}%
\]
where $r=2k$ or $2k+1$ for $\mathbf{2}$\textbf{-}$\mathbf{D}_{\mathbf{2k}}$
and $\mathbf{1}$\textbf{-}$\mathbf{D}_{\mathbf{2k+1}}^{\mathbf{II}}$
respectivly. We claim that $u_{0}^{s-1}v_{0}^{s},$ $s\geq0,$ generates
$H^{1}(\widetilde{U}_{RDP,i},\mathcal{O}_{Z-E}(2Z))\cong H_{k-1}.$ This can be
seen for instance by proving that the class of the element $(v_{0}+1)^{t}%
u_{0}^{2t-1}$ equals the class of $u_{0}^{2t-1}v_{0}^{2t}$ in $H_{t+1}$ and
that this element is nonzero and generates the kernel of $H_{t+1}\rightarrow
H_{t}.$ The claim follows by induction. We have $u_{0}^{s-1}v_{0}^{s}%
=z_{3}^{s+1}/z_{1}.$ Since $1/v_{0}\in\mathbb{C}[u_{0},v_{0},u_{0}%
^{-1}]/((v_{0}+1)^{t+1})$, we will have $1/z_{1}\in B^{\prime}/(x)$, so in the
$\mathbf{2}$\textbf{-}$\mathbf{D}_{\mathbf{2k}}$-case we may define
$D\in\operatorname*{Der}(P,B^{\prime})$ by $D(z_{2})=$ $z_{3}^{s}/z_{1}$,
$D(z_{1})=D(z_{3})=0.$ This gives $D(x)=z_{3}^{s+1}/z_{1}$ and $D(g)=z_{3}%
^{s}.$ In the $\mathbf{1}$\textbf{-}$\mathbf{D}_{\mathbf{2k+1}}^{\mathbf{II}}%
$-case we define $D$ by $D(z_{2})=$ $z_{3}^{s+1}/z_{1}$, $D(z_{1}%
)=D(z_{3})=0,$ and we get $D(x)=D(z_{2})=z_{3}^{s+1}/z_{1}$ and $D(g)=z_{3}%
^{s+1}.$ We thus get $\operatorname*{coker}\alpha_{i}=0$ in both cases.

\textbf{Remaining cases: }The remaining cases are 1-configurations, and may be
checked in similar fashion. However, these cases also follow from \cite{DJ},
so they are omitted.
\end{proof}

We may now prove the main theorem stated in the introduction. \ In fact, we
state and prove a slightly more general version:

\begin{theorem}
\label{etteo}Let $X$ be a rational surface singularity with embedding
dimension $e\geq4$ and with almost reduced fundamental cycle. Then $c(X)$ is
less or equal to the number of $\mathbf{3}$\textbf{-}$\mathbf{A}_{\mathbf{n}}%
$-configurations in the dual graph and greater or equal to the number of
$\mathbf{3}$\textbf{-}$\mathbf{A}_{\mathbf{n}}$-configurations with the
property thet the adjacent non $-2$-curves intersects $Z$ negatively.
\end{theorem}

\begin{proof}
The possible dual graphs for $X$ are classified in proposition \ref{badconfig}%
. From lemma \ref{alphalemma} there is a surjection $\mathbb{C}^{s}\rightarrow
H^{1}(\mathfrak{m}\mathcal{F}^{1})$ where $s$ is the number $\mathbf{3}%
$\textbf{-}$\mathbf{A}_{\mathbf{n}}$-configurations in the dual graph. It
follows from the proof of lemma \ref{alphalemma} that the map is injective
when we restrict to the copies of $\mathbb{C}$ corresponding to
RDP-configurations which blow down to RDPs. If the fundamental cycle $Z$
intersects the three non $-2$-curves adjacent to a $\mathbf{3}$\textbf{-}%
$\mathbf{A}_{\mathbf{n}}$-configuration negatively, this will be the case.
\end{proof}

\begin{corollary}
Let $X$ be a rational determinantal or quasi-determinantal surface singularity
with embedding dimension $e\geq4.$ Then $c(X)=0.$
\end{corollary}

\begin{remark}
Theorem \ref{etteo}, generalizes and \textit{corrects} Theorem 2.6.3 of
\cite{mythesis} where some $2$-configurations mistakenly are computed to
contribute positively to $c(X).$
\end{remark}

\providecommand{\bysame}{\leavevmode\hbox to3em{\hrulefill}\thinspace}
\providecommand{\MR}{\relax\ifhmode\unskip\space\fi MR }
\providecommand{\MRhref}[2]{%
  \href{http://www.ams.org/mathscinet-getitem?mr=#1}{#2}
}
\providecommand{\href}[2]{#2}

\end{document}